\documentclass[12pt]{article}
\usepackage[cp1251]{inputenc}
\usepackage{amssymb}                         
\usepackage{color}                           
\usepackage[unicode]{hyperref}               
\ifx\pdfoutput\undefined                     
\usepackage{graphicx}
 \else
\usepackage[pdftex]{graphicx}
\usepackage{cite}
\fi
\usepackage[russian]{babel}
\begin{document}

\pagestyle{plain}

\begin{center}
{\LARGE Модель с непрерывным представлением времени для\\
\par\vspace{0.2ex} задачи
составления расписаний многопродуктового\\
\par\vspace{1.2ex} производства}\\
\mbox{ }\\
\mbox{ }\\
{\large Ю.В.~Коваленко}\\
\mbox{ }\\
{\small Омский государственный университет им. Ф.М. Достоевского,
пр. Мира, 55а, 644077}\\
{\small Омск, Россия, E-mail: juliakoval86@mail.ru}\\
\mbox{ }\\
\end{center}

%
\mbox{ }\\
{\small \textbf{Аннотация.} Рассматривается задача составления
расписаний многопродуктового производства. Особенностью постановки
является то, что каждый продукт имеет несколько технологий
производства, при выполнении которых используется сразу несколько
машин, работающих одновременно. При этом, если машина
переключается с одной технологии на другую, то необходимо
выполнять переналадку. Построены модели частично целочисленного
линейного программирования для задачи в общей постановке и для
случая, когда длительности переналадки удовлетворяют неравенству
треугольника. Для сравнения предложенных моделей проведены
численные
эксперименты на построенных случайным образом тестовых примерах.}\\
\\
{\small \textbf{Ключевые слова:} расписание, переналадки,
технологии, целочисленное линейное программирование.}\\

{

\section{Введение}
\label{sec:problemdefinition}
\parindent=1cm
\vspace{-0.5cm} \hspace{0cm}
\par
Рассмотрим задачу составления расписаний многопродуктового
производства (СРМП) следующего вида. Имеется химическое
предприятие, выпускающее $k$ различных продуктов. ${{V_{i}}\in
{\mathbb R}^+}$~--~требуемый объем производства продукта $i,\
{i=1,\dots,k}$. Здесь и далее ${\mathbb R}^+$ обозначает множество
положительных вещественных чисел. Пусть $m$~--~число машин,
которые могут использоваться при выпуске продукции.

Для каждого продукта $i,\ i=1,\dots,k,$ указана одна или более
технологий его производства. Пусть $U$~--~множество технологий,
где $|U|=d$, каждая из которых характеризуется набором
одновременно занимаемых машин ${M_u\subseteq \{1,\dots,m\}, \
{u\in U}},$ т.~е. если производится продукт $i$ по технологии $u$,
то одновременно задействованы все машины, относящиеся к данной
технологии. При этом в любой момент  каждая машина не может быть
задействована  более чем в одной технологии. Прерывание выполнения
технологий допускается.

Пусть $U_i \subseteq U$~--~множество технологий по производству
продукта $i,$ ${i=1,\dots, k}$, для каждой из которых  задан объем
$a_{u}\in {\mathbb {R}}^+$ выпуска данного продукта в единицу
времени, ${u\in U_{i}}.$ Предполагается, что для выпуска
продукта~$i$ может быть использовано несколько технологий из
множества~$U_i$, $i=1,\dots,k$.

Для машины $l$ заданы длительности переналадки этой машины с
технологии~$u$ на технологию~$q$, обозначенные через $s_{l u q}
\in {\mathbb R}_+$  (здесь  ${\mathbb R}_+$ -- множество
неотрицательных вещественных чисел), для всех ${ u,\ q\in K_{l},}$
где ${K_l=\{u:\ l\in M_u,\ u \in U\}}$ -- множество технологий,
использующих машину~$l$$,\ {l=1,\dots,m}$.

Для каждого продукта $i,\ i=1,\dots,k,$ необходимо определить,
какие технологии~${u\in U_{i}}$ будут использоваться для его
производства, и для выбранных технологий составить расписание их
выполнения с учетом переналадок машин и невозможности
одновременного использования одной машины в различных технологиях,
таким образом, чтобы общее время окончания производства всех
продуктов~$C_{\max}$ в объемах~$V_1,\dots,V_k $ было минимально.

Если длительности переналадки удовлетворяют следующему
неравенству:

\begin{equation}\label{triengle} s_{luq}+s_{lqp}\ge s_{lup},\
l=1,\dots,m,\  u,q,p\in K_l,\end{equation} то будем говорить, что
они удовлетворяют  неравенству треугольника. Отметим, что на
практике достаточно часто встречаются задачи составления
производственных расписаний, в которых длительности переналадки
удовлетворяют неравенству треугольника.

Поставленная задача является NP-трудной в сильном смысле, так как
в частном случае при ${m=1}$ и длительностях переналадки,
удовлетворяющих  неравенству треугольника, к ней сводится
полуметрическая задача о кратчайшем гамильтоновом пути, являющаяся
NP-трудной в сильном смысле~\cite{IPS}.

В настоящей статье построена модель частично целочисленного
линейного программирования для задачи СРМП в общей постановке и
для случая, когда длительности переналадки удовлетворяют
неравенству треугольника. Вторая модель представляется более
предпочтительной при использовании методов оптимизации, основанных
на ЛП-релаксации, что подтверждается численным экспериментом на
сгенерированных случайным образом тестовых примерах.

\section{Модель частично целочисленного линейного
программирования}

\label{sec:problemmodel}
\parindent=1cm
\vspace{-0.5cm} \hspace{0cm}
\par
Известно множество подходов к формулированию задач построения
производственных расписаний в виде моделей частично целочисленного
линейного программирования, с учетом сложных связей между
технологиями, продуктами и машинами~\cite{BPA, FKPS, IF}. Учитывая
структуру рассматриваемой задачи, можно предложить модель,
основанную на тех же принципах.

Перед описанием модели определим понятие точки событий,
аналогичное введенному  в работе~\cite{IF}. Точка событий~--~это
группа переменных  задачи, которые задают некоторый набор
технологий и моменты начала и окончания выполнения технологий из
этого набора. При этом в одной точке событий каждая машина может
быть задействована не более чем в одной технологии.


Введем следующие обозначения:\\
$I$~--~множество продуктов$,\ |I|=k;$\\
$L$~--~множество машин$,\ |L|=m;$\\
$N=\{1,\dots, n_{\max}\}$~--~множество
точек событий;\\
$H= \sum\limits_{i\in I}\max\limits_{u\in
U_{i}}\left\{\frac{V_i}{a_u}\right\}+(k-1)\cdot \max\limits_{l\in
L,\ u,q\in K_l}\{s_{l u q}\}$~--~оценка сверху длины расписания$.$
Ясно, что данного времени  достаточно для производства всех
продуктов.

Определим переменные задачи:\\
$w_{u n}$~--~бинарная переменная, которая определяет, выполняется
технология~$u\in U$ в точке событий~$n\in N$ (при этом $w_{u
n}=1$) или нет
($w_{u n}=0$);\\
 $y_{l n}$~--~бинарная переменная, которая определяет, используется машина~$l\in L$ в точке событий~$n\in N$ (при этом $y_{l n}=1$) или нет
($y_{l n}=0$);\\
$T^s_{u  n}$~--~время начала выполнения технологии~$u\in U$ в
точке
событий~$n\in N$;\\
$T^f_{u n}$~--~время завершения выполнения технологии~$u\in U$ в
точке
событий~$n\in N$;\\
$C_{\max}$~--~момент завершения производства всех продуктов$.$

С использованием введенных обозначений и переменных модель
частично целочисленного линейного программирования для задачи СРМП
может быть записана следующим образом:

\begin{equation}\label{1''N} C_{\max}\to \min, \end{equation}

\begin{equation}
 \label{2N}
T^f_{u  n}\leqslant C_{\max}, \ u\in U,\
 n \in N,
 \end{equation}

\begin{equation}
 \label{3N}
\sum_{u \in K_l}w_{u n}=y_{l n},\  l\in L,\  n\in N,
 \end{equation}


 $$T^s_{u  n}\geqslant
T^f_{q  \tilde{n}}+s_{l q u}-H\cdot(2-w_{u n}-w_{q
\tilde{n}}+\sum_{\tilde{n}<n'<n}{y_{l n'}}),$$
\begin{equation}
 \label{6N}
  l\in L,\
u,q\in K_l,\  n,\ \tilde{n} \in N,\ n\ne 1, \ \tilde{n}<n,
 \end{equation}

\begin{equation}
 \label{9N}
T^f_{u  n} \geqslant T^s_{u  n},\   u\in U, \  n\in N,
 \end{equation}

 \begin{equation}
 \label{11N}
T^f_{u  n} - T^s_{u  n} \leqslant w_{u n}\cdot \max_{q\in
U_i}\left\{\frac{V_i}{a_{q}}\right\} ,
 \ i\in I,\  u\in U_i, \  n\in N,
 \end{equation}

\begin{equation}
 \label{12N}
\sum_{n\in N}\sum_{u\in U_i}a_u\cdot (T^f_{u  n} - T^s_{u  n})
\geqslant V_i ,
 \ i\in I,
 \end{equation}

  \begin{equation}
 \label{10N}
  T^s_{u  n} \geqslant 0,\ u\in U, \    n\in N,
 \end{equation}

 \begin{equation}
 \label{15N}
w_{u n }\in\{0,1\} , \ u\in U,\  n \in N,
 \end{equation}

 \begin{equation}
 \label{15'N}
y_{l n }\in\{0,1\} , \ l\in L,\  n \in N.
 \end{equation}

{Целевая функция~(\ref{1''N}) и неравенство~(\ref{2N})} задают
критерий  минимизации момента окончания производства всех
продуктов. Ограничение~(\ref{3N}) выражает то, что в каждой точке
событий {машина~$l$}  либо используется, причем не более чем в
одной технологии, если ${y_{l n}=1}$, либо нет в противном случае.
Ограничение~(\ref{6N}) говорит о том, что время начала
технологии~$u$ на машине~$l$ должно быть не меньше, чем время
окончания предыдущей технологии на той же машине плюс длительность
переналадки. Условие~(\ref{9N}) гарантирует неотрицательность
длительности технологий. Если технология~$u$ в точке событий~$n$
не выполняется, т.~е. ${w_{u n}=0}$, то ее длительность должна
быть равна нулю, что обеспечивается {условием~(\ref{11N}).}
{Ограничение~(\ref{12N})} гарантирует выпуск продукции в заданном
объеме. Ограничения~(\ref{10N})~--~(\ref{15'N}) описывают область
определения переменных.

На практике достаточно часто возникают задачи составления
производственных расписаний, в которых длительности переналадки
удовлетворяют неравенству треугольника~(\ref{triengle}). Поэтому
имеет смысл сформулировать  модель частично целочисленного
линейного программирования для задачи СРМП при выполнении
неравенства треугольника. Предлагаемая модель представляется более
предпочтительной, так как в ней удается исключить
неравенство~(\ref{6N}), которое усложняет задачу для методов
оптимизации, основанных на ЛП-релаксации.

Используя введенные ранее обозначения,
построим модель частично целочисленного линейного программирования в следующем виде:\\

\begin{equation}\label{1''} C_{\max}\to \min, \end{equation}

\begin{equation}
 \label{2}
T^f_{u  n}\leqslant C_{\max}, \ u\in U,\
 n \in N,
 \end{equation}

\begin{equation}
 \label{3}
\sum_{u \in K_l}w_{u n}=y_{ln},\  l\in L,\  n\in N,
 \end{equation}

\begin{equation}
 \label{4'}
T^s_{u  ,n+1}\geqslant T^f_{u  n}, \ u\in U,\  n \in N,\ n\ne
n_{\max},
 \end{equation}

 $$T^s_{u  ,n+1}\geqslant
T^f_{q  n}+s_{l q u}\cdot w_{u ,n+1}-H\cdot(1-w_{u ,n+1}),$$
\begin{equation}
 \label{6}
  l\in L,\
u,q\in K_l,\ u \neq q,\  n \in N,\ n\ne n_{\max},
 \end{equation}

  \begin{equation}
 \label{10}
  T^s_{u  n} \geqslant -H\cdot(1-w_{u n}),\ u\in U, \    n\in N,
 \end{equation}

\begin{equation}
 \label{9}
T^f_{u  n} \geqslant T^s_{u  n},\   u\in U, \  n\in N,
 \end{equation}

 \begin{equation}
 \label{11}
T^f_{u  n} - T^s_{u  n} \leqslant w_{u n}\cdot \max_{q\in
U_i}\left\{\frac{V_i}{a_{q}}\right\} ,
 \ i\in I,\  u\in U_i, \  n\in N,
 \end{equation}

\begin{equation}
 \label{12}
\sum_{n\in N}\sum_{u\in U_i}a_u\cdot (T^f_{u  n} - T^s_{u  n})
\geqslant V_i ,
 \ i\in I,
 \end{equation}

 \begin{equation}
 \label{15'}
w_{u n }\in\{0,1\} , \ u\in U,\  n \in N,
 \end{equation}

 \begin{equation}
 \label{15}
y_{l n }\in\{0,1\} , \ l\in L,\  n \in N.
 \end{equation}

{Ограничения~(\ref{1''})~--~(\ref{3}) и  (\ref{9})~--~(\ref{15})
имеют тот же смысл, что и в предыдущей модели}.
Ограничение~(\ref{4'}) выражает то, что время начала выполнения
технологии~$u$ в точке событий~$n+1$ должно быть не меньше, чем
время окончания ее выполнения в точке событий~$n$.
Ограничение~(\ref{6}) говорит о том, что время начала
технологии~$u$ на машине~$l$ должно быть не меньше, чем время
окончания предыдущей технологии на той же машине плюс длительность
переналадки.  Однако, если технология~$u$ является первой на
машине~$l$, то длительность переналадки на нее учитывать не нужно.
Это обеспечивается благодаря тому, что переменные~$T^s_{q n}$  во
всех предшествующих точках событий могут принимать отрицательные
значения. Если же технология $u$ имеет место в точке событий~$n$,
т.~е. ${w_{u n}=1}$, то время начала ее использования должно быть
неотрицательным, что обеспечивается условием~(\ref{10}).

Необходимо отметить, что  условие (\ref{6}) будет гарантировать
получение оптимального решения задачи, только если длительности
переналадки удовлетворяют неравенству треугольника.

\section{Вычислительный эксперимент}
\label{sec:experiment}
\parindent=1cm
\vspace{-0.5cm}
 \hspace{0cm}
\par
Для сравнения предложенных моделей частично целочисленного
линейного программирования был проведен вычислительный эксперимент
на случайным образом сгенерированных задачах трех серий $S1,\ S2$
и $S3$.

Модели (\ref{1''N})~--~(\ref{15'N}) и (\ref{1''})~--~(\ref{15})
были записаны в системе моделирования GAMS~22.6 \ и задача СРМП
решалась с помощью универсального пакета~CPLEX~11.0. При этом
использовался метод ветвей и границ с отсечениями, настройки
которого были выбраны по умолчанию. Тестирование проводилось на
ЭВМ Intel~Core2~Duo~CPU~E7200~2.54~ГГц, оперативная память~2 Гб.

При генерации тестовых задач  для каждой технологии ${u\in U}$
число машин~${|M_u|\in \{1,\dots,m\}}$ выбиралось случайно с
равномерным распределением, а затем машины случайно назначались на
данную технологию без повторений. Числовые значения для всех
тестовых задач генерировались случайным образом с равномерным
распределением из следующих множеств: ${|U_i|\in \{1,\dots,
U_{\max}\}}$; $V_i\in [1,V_{\max}]$; $a_u\in [1,\frac{V_i}{2}]$,
где $i$ такое, что~$u\in U_i$; ${s_{luq}\in [0,s_{\max}]}$. В
табл.~\ref{tabser} приведены выбранные значения параметров для
каждой серии.

\begin{table}[!h]
\begin{center}
\caption{}
\small{Параметры серий\\} \vspace{0.3cm} \label{tabser}
\begin{tabular}{|c|c| c| c| c|c| c| c|c|}
\hline
{серия} &{число задач} & $k$ & $m$ & $U_{\max}$ & $V_{\max}$ & $s_{\max}$  & $ n_{\max}$ \\
\hline
$S1$ & $10$ & $4$ & $4$ & $3$ & $10$ & $5$ &  $4$ \\
\hline
$S2$ & $10$ & $5$ & $7$ & $5$ & $12$ & $7$ &  $5$ \\
\hline
$S3$ & $10$ & $7$ & $9$ & $6$ & $15$ & $9$ &  $7$ \\
\hline
\end{tabular}
\end{center}
\end{table}

Для проведения вычислительного эксперимента было установлено
максимальное время равное 5000~сек., которое CPLEX может решать
одну тестовую задачу.

При описании результатов вычислительного эксперимента будут использоваться следующие обозначения:\\
$VAR$~--~количество переменных в моделях (\ref{1''N})~--~(\ref{15'N}) и (\ref{1''})~--~(\ref{15});\\
$EQV$~--~число ограничений в модели (\ref{1''N})~--~(\ref{15'N});\\
$EQV_{\Delta}$~--~число ограничений в модели (\ref{1''})~--~(\ref{15});\\
$C$~--~значение целевой функции, полученное при
решении CPLEX  задачи~(\ref{1''N})~--~(\ref{15'N});\\
$C_{\Delta}$~--~значение целевой функции, полученное при
решении CPLEX  задачи~(\ref{1''})~--~(\ref{15});\\
$t$~--~время работы CPLEX в сек. при
решении задачи (\ref{1''N})~--~(\ref{15'N});\\
$t_{\Delta}$~--~время работы CPLEX в сек. при решении задачи
(\ref{1''})~--~(\ref{15}).

\begin{table}[!h]
\begin{center}
\caption{}
\small{Сравнение моделей на задачах серии~$S1$\\}
\vspace{0.3cm} \label{tabms1}
\begin{tabular}{|c|c| c| c| c|c| c| c|}
\hline
{\bf задача} & $VAR$ & $EQV$ & $EQV_{\Delta}$ & $C$  & $C_{\Delta}$ & $t$ & $t_{\Delta}$\\
\hline
 1 &  113 & 530& 334 & 19.398 & 19.398 &  4.3    & 2.7\\
\hline
 2 & 137 & 896 & 528 &  19.628 & 19.628 &   3 &  1.8\\
\hline
 3 &  113 & 662 & 388 & 21.757 & 21.757  &  2.4   & 1\\
 \hline
 4 &  101 & 380 & 255 & 14.401 & 14.401 & 3.6  &  1.2\\
 \hline
 5 &  113  & 584 & 358 & 19.678 & 19.678 & 4.8 & 2.7\\
 \hline
 6 & 89 & 290 & 200  & 14.958  & 14.958  &  2.2 & 1.1\\
 \hline
 7 & 101 & 242 & 195 & 9.675 & 9.675   & 3.1  & 2.2\\
 \hline
 8 & 113 & 614 & 376 & 25.583 & 25.583 & 4.4 & 2.4\\
 \hline
 9 & 101 & 416 & 267 & 15.309 & 15.309 & 5.5 & 2.4\\
 \hline
 10 & 113 & 542 & 340 & 22.818 & 22.818 & 4.6 & 1.2\\
 \hline
\end{tabular}
\end{center}
\end{table}

\begin{table}[!h]
\begin{center}
\caption{}
\small{Сравнение моделей на задачах серии~$S2$\\}
\vspace{0.3cm} \label{tabms2}
\begin{tabular}{|c|c| c| c| c|c| c| c|}
\hline
{\bf задача} & $VAR$ & $EQV$ & $EQV_{\Delta}$ & $C$  & $C_{\Delta}$ & $t$ & $t_{\Delta}$\\
\hline
1 &  231 & 2485 & 1112 & 24.274 & 24.274 &  143 & 51\\
\hline
2 &  231 & 3075 & 1312 & 33.079 & 33.079 &  23  & 18\\
\hline
3 &  246 & 3820 & 1616 & 41.768 & 41.768 &  113 & 24\\
\hline
4 & 306 & 5280 & 2232 & 14.505 & 14.505 & 331 & 111\\
\hline
5 & 231 & 3675 & 1536 & 37.327 & 37.327 & 1112 & 655\\
\hline
6 & 321 & 6385 & 2672 & 20.463 & 20.463 & 2000  & 486\\
\hline
7 & 321 & 5005 & 2152 & 16.211 & 16.211 & 1013 & 216\\
\hline
8 & 261 & 2665 & 1200 & 13.437 & 13.437 & 978 & 419\\
\hline
9 & 231 & 2635 & 1160 & 32.599 & 32.599 & 526  & 212\\
\hline
10 & 201 & 1945 & 864 &46.936 & 46.936 & 108 & 45\\
 \hline
\end{tabular}
\end{center}
\end{table}

\begin{table}[!h]
\begin{center}
\caption{}
\small{Сравнение моделей на задачах серии~$S3$\\}
\vspace{0.3cm} \label{tabms3}
\begin{tabular}{|c|c| c| c| c|c| c| c|}
\hline
{\bf задача} & $VAR$ & $EQV$ & $EQV_{\Delta}$ & $C$  & $C_{\Delta}$ & $t$ & $t_{\Delta}$\\
\hline
 1 &  505 & 15862 &   4720 &     29.133 & 29.133 & 2560 & 1000\\
\hline
 2 &  505 & 17143 &  5056  &   57.523 &  57.523 & 1360 & 700\\
\hline
 3 &  568 & 23044 &  6742  &   38.543 & 38.543 & 3053 & 1987\\
\hline
 4 &  526 & 15085 &  4550  &   37.267 &  37.267 & 5000 & 2876\\
\hline
 5 &  631 &  22414 &  6652  &   25.419 &  25.419 & 5000 & 1159\\
\hline
 6 &  568 & 22057 &  6466  &   36.557 & 36.557 & 5000 & 1897\\
\hline
 7 &  694 & 29911 &  8782  &   {48.996$^*$} & { 43.703} & 5000 & 3001\\
\hline
 8 &  568 & 22960 &  6718  &   43.610 & 43.610 & 5000 & 2564\\
\hline
 9 &  526 & 18088 &  5330  &   {46.412$^*$} & { 46.342} & 5000 & 2390\\
\hline
 10 & 631 & 31507 &  9148  &   {59.180$^*$} & { 58.539} & 5000 &  1602\\
\hline

\end{tabular}
\end{center}
\end{table}

В табл. \ref{tabms1}, \ref{tabms2} и \ref{tabms3} представлены
результаты вычислительного эксперимента. Видно, что модель
(\ref{1''})~--~(\ref{15}) содержит меньшее число ограничений, чем
модель (\ref{1''N})~--~(\ref{15'N}). Кроме того, в среднем на
сериях $S1$ и $S2$ пакет CPLEX решает тестовые примеры, записанные
в модели (\ref{1''})~--~(\ref{15}), более чем  в два раза быстрее,
чем при их записи в модели (\ref{1''N})~--~(\ref{15'N}). На
серии~$S3$ при решении пакетом CPLEX задач 4~--~10, записанных в
модели (\ref{1''N})~--~(\ref{15'N}),  за установленное время
удается найти только допустимые решения, которые не всегда
являются оптимальными (отмечены '$*$').

Таким образом, если в задаче СРМП длительности переналадки
удовлетворяют неравенству треугольника, то при ее решении лучше
использовать модель (\ref{1''})~--~(\ref{15}), так как данная
модель имеет меньшее число ограничений и является менее сложной
для методов оптимизации, основанных на ЛП - релаксации (в том
числе для пакета CPLEX).
\section{Заключение}\label{sec:concl}
\parindent=1cm
\vspace{-0.5cm}
 \hspace{0cm}
\par
 В работе рассмотрена задача составления расписаний
 многопродуктового производства. Особенностью данной задачи является то, что каждый
 продукт может производиться по нескольким технологиям, каждая из
 которых характеризуется набором одновременно занимаемых машин.
Построены модели частично целочисленного линейного
программирования для задачи в общей постановке и для случая, когда
длительности переналадки удовлетворяют неравенству треугольника. С
помощью вычислительного эксперимента показано, что вторая модель
является более предпочтительной для методов оптимизации,
основанных на ЛП-релаксации.\\

Автор благодарит А.В.~Еремеева за предложенную постановку задачи.

\renewcommand{\refname}{ЛИТЕРАТУРА}



\begin {thebibliography}{}

\bibitem{BPA} \textbf{Борисовский П.А.} {Генетический алгоритм для одной задачи
составления производственного расписания с переналадками}~// Тр.
XIV Байкальской международной школы-семинара <<Методы оптимизации
и их приложения>>. -- Иркутск:  ИСЭМ СО РАН, 2008. Т.~4.
С.~166~--~173.

\bibitem{FKPS}  \textbf{Floudas C.A., Kallrath J., Pitz H.J., Shaik M.A.}
{Production scheduling of a large-scale industrial continuous
plant: short-term and medium-term scheduling}~//
Comp.~Chem.~Engng. 2009. Vol.~33. P.~670~--~686.

\bibitem{IF} \textbf{Ierapetritou M.G, Floudas C.A.}  {Effective
continuous-time formulation for short-term scheduling: I.
multipurpose batch process}~// Ind.~Eng.~Chem.~Res.
 1998.  Vol.~37.  P.~4341~--~4359.
\bibitem{IPS} \textbf{Itai A., Papadimitriou C.H., Szwarcfiter
J.L.} {Hamilton paths in grid graphs}~// SIAM Journal of
Computing.  1982.  Vol.~11, N 4. P.~676~--~686.


\end{thebibliography}
\end{document}